\documentclass{article}
\usepackage[utf8]{inputenc}
\usepackage{amssymb,amsmath}
\usepackage[usenames, dvipsnames]{color}
\usepackage{amsthm}
\usepackage{centernot}
\usepackage{mathtools}
\usepackage{stmaryrd}
\usepackage{array}
\usepackage{float} 
\usepackage{setspace}
\usepackage{longtable}
\newtheorem{theorem}{Theorem}

\newtheorem{lemma}{Lemma}

\begin{document}
\begin{center}
{\Large Classification of $\lambda$-homomorphic braces on $\mathbb{Z}^2$}

\bigskip T.~Nasybullov, I. Novikov
\end{center}
\begin{abstract}
If $A=(A,\oplus,\odot)$ is a $\lambda$-homomorphic brace with $(A,\oplus)=\mathbb{Z}^2$, then the operations in this brace are given by formulas
\begin{align*}\begin{pmatrix}a_1\\a_2\end{pmatrix}\oplus\begin{pmatrix}b_1\\b_2\end{pmatrix}=\begin{pmatrix}a_1+b_1\\a_2+b_2\end{pmatrix},&&\begin{pmatrix}a_1\\a_2\end{pmatrix}\odot\begin{pmatrix}b_1\\b_2\end{pmatrix}=\begin{pmatrix}a_1\\a_2\end{pmatrix}+\varphi^{a_1}\psi^{a_2}\begin{pmatrix}b_1\\b_2\end{pmatrix},
\end{align*}
where $\varphi,\psi\in{\rm GL}_2(\mathbb{Z})$ are cpecific matrices which depend on $A$. Not every pair $(\varphi,\psi)$ lead to a brace. In the present paper we find all possible pairs $(\varphi,\psi)$ of matrices from  ${\rm GL}_2(\mathbb{Z})$ which lead to $\lambda$-homomorphic braces with $(A,\oplus)=\mathbb{Z}^2$. The obtained result gives the full classification of $\lambda$-homomorphic braces on $\mathbb{Z}^2$ which was started by Bardakov, Neshchadim and Yadav in \cite{BarNesYad}.

 ~
 
 ~
\end{abstract}

\section{Introduction}
A set-theoretical solution of the Yang-Baxter equation on a set $X$ is a map $r:X\times X\to X\times X$ which satisfies the equality 
$$(r\times id)(id \times r)(r\times id)=(id\times r)(r \times id)(id\times r)$$
called the set-theoretical Yang-Baxter equation. This equation comes from the quantum Yang-Baxter equation which appeared in theoretical physics and statistical mechanics in the works of Yang \cite{Yan} and Baxter \cite{Bax1,Bax2}. The set-theoretical Yang-Baxter equation has various applications in different fields of mathematics. For example, solutions of this equation can be used for constructing representations of braid-like groups and invariants of classical and virtual knot theories \cite{BarNasmult1, BarNasmult2, BarNasmult3}. 

The problem of classifying set-theoretical solutions of the Yang-Baxter equation was formulated by Drinfel'd in \cite{Dri}. One approach to this problem is to classify solutions of the Yang-Baxter equation in the case when $X$ is not just a set, but some algebraic system. Many new algebraic systems have been built along this way. Examples of such algebraic systems are braces introduced by Rump in \cite{Rum} and skew braces introduced by Guarnieri and Vendramin in \cite{GuaVen}. 

A skew brace $A=(A,\oplus,\odot)$ is an algebraic system with two binary algebraic operations  $\oplus$, $\odot$ such that $A_{\oplus}=(A,\oplus)$, $A_{\odot}=(A,\odot)$ are groups and the equality
\begin{equation}\label{mainleft}
a\odot(b\oplus c)=(a\odot b)\oplus(\ominus a) \oplus (a\odot c)
\end{equation}
holds for all $a,b,c\in A$, where $\ominus a$ denotes the inverse to $a$ element with respect to the operation~$\oplus$. The group $A_{\oplus}$ is called the additive group of a skew brace $A$, and the group $A_{\odot}$ is called the multiplicative group of a skew brace $A$. If $A_{\oplus}$ is abelian, then $A$ is called a brace. Skew braces have connections with other algebraic structures such as groups with exact factorizations, Zappa-Sz\'{e}p products, triply factorized groups and Hopf-Galois extensions  \cite{SmoVen}. Some algebraic aspects of skew braces are studied in \cite{CedSmoVen, Chi, GorNas, Jes, KonSmoVen, Nas, SmoVen}. A big list of problems concerning skew braces is collected in \cite{Ven}.

If $A=(A,\oplus,\odot)$ is a skew brace, then the map $r:A\times A\to A\times A$ given by 
\begin{equation}\label{rumybe}
r(x,y)=\Big(\ominus x\oplus(x\odot y),(\ominus x\oplus(x\odot y))^{-1}\odot x\odot y\Big)
\end{equation}
for $x,y\in A$ is a set-theoretical solution of the Yang-Baxter equation on $A$. If in this situation $A$ is a brace, then the solution $r$ is non-degenerate and involutive. Moreover, every non-degenerate involutive solution of the Yang-Baxter equation can be given by formula (\ref{rumybe}) for a suitable brace~$A$ \cite{Rum}. Thus the question of classification of all set-theoretical solutions of the Yang-Baxter can be reformulated in a weaker form as the question of classification of all skew braces.

At the moment, there are quite a few classification results about finite skew braces. In \cite{GuaVen} the algorithm of enumerating and constructing finite skew braces is built. This algorithm is used to produce a database of skew braces of small size. In \cite{BarNesYad0} this algorithm is modified, and the modification is used for enumerating finite skew braces of orders up to $868$ with some exceptions. In \cite{bonclass0} all skew braces whose orders are products of two distinct primes are described. In \cite{alibyot} all skew braces whose orders are products of three distinct primes are classified. Papers \cite{bonclass1, bonclass2, Campedel, Campedel2} give the classification of skew braces of order $p^2q$, where $p, q$ are distinct prime numbers. Paper \cite{BEJP} gives the  description of all possible multiplicative groups of finite skew braces whose additive group has trivial centre. 

Far fewer classification results are known for infinite skew braces. Little is known even in the case when the additive group $A_{\oplus}$ is the free abelian group $\mathbb{Z}^n$. In \cite{Rum2} it is proved that if  $(A, \oplus, \odot)$ is a brace with $(A, \oplus) = \mathbb{Z}$, then the operation $\odot$ is given by one of the following formulas
\begin{enumerate}
\item $a \odot b = a + b$ for all $a,b\in A$ (in this case $A_{\odot}=A_{\oplus}=\mathbb{Z}$),
\item $a \odot b = a + (-1)^{a}b$ for all $a,b\in A$ (in this case $A_{\odot}=\mathbb{Z}\rtimes \mathbb{Z}_2$).
\end{enumerate}
In \cite{BarNesYad} some examples of $\lambda$-homomorphic skew braces $(A,\oplus,\odot)$ with $A_{\oplus}=\mathbb{Z}^2$ are found. More precisely, it is proved that if 
\begin{align}\label{bardakovresultsyouknow}
\varphi=\psi=\begin{pmatrix}1+m&m\\-m&1-m\end{pmatrix}&&\text{or}&&\varphi=\psi=\begin{pmatrix}1+m&2+m\\-m&-1-m\end{pmatrix}
\end{align}
for some $m\in \mathbb{Z}$, then the algebraic system $(\mathbb{Z}^2,\oplus,\odot)$ with 
\begin{align*}\begin{pmatrix}a_1\\a_2\end{pmatrix}\oplus\begin{pmatrix}b_1\\b_2\end{pmatrix}=\begin{pmatrix}a_1+b_1\\a_2+b_2\end{pmatrix},&&\begin{pmatrix}a_1\\a_2\end{pmatrix}\odot\begin{pmatrix}b_1\\b_2\end{pmatrix}=\begin{pmatrix}a_1\\a_2\end{pmatrix}+\varphi^{a_1}\psi^{a_2}\begin{pmatrix}b_1\\b_2\end{pmatrix},
\end{align*}
for all $a_1,a_2,b_1,b_2\in \mathbb{Z}$ is a $\lambda$-homomorphic brace with $A_{\oplus}=\mathbb{Z}^2$. 

In the present paper we complete the classification of $\lambda$-homomorphic braces $(A,\oplus,\odot)$ with $A_{\oplus}=\mathbb{Z}^2$ started in \cite{BarNesYad}. The main result of the paper is the following theorem. The denotation 
$A\equiv_kB$ used in the formulation of the theorem means that $A=B+kC$ for appropriate matrix $C$.\medskip

\noindent \textbf{Theorem A.} {\it Let $(A,\oplus,\odot)$ be a $\lambda$-homomorphic brace with $A_{\oplus}=\mathbb{Z}^2$. Then the operations in this brace are given by formulas
\begin{align*}\begin{pmatrix}a_1\\a_2\end{pmatrix}\oplus\begin{pmatrix}b_1\\b_2\end{pmatrix}=\begin{pmatrix}a_1+b_1\\a_2+b_2\end{pmatrix},&&\begin{pmatrix}a_1\\a_2\end{pmatrix}\odot\begin{pmatrix}b_1\\b_2\end{pmatrix}=\begin{pmatrix}a_1\\a_2\end{pmatrix}+\varphi^{a_1}\psi^{a_2}\begin{pmatrix}b_1\\b_2\end{pmatrix},
\end{align*}
where $\varphi,\psi\in{\rm GL}_2(\mathbb{Z})$ are matrices which satisfy one of the conditions listed in Table 1.}

{\footnotesize
\renewcommand{\arraystretch}{1.5}
\begin{longtable}{|p{1em}|p{5.5em}|p{2em}|p{20em}|p{15em}|}
\hline
~&~&~ &$\varphi$&$\psi$\\\cline{1-5}
1.&${\rm det}(\varphi)=1$&1.1&$\varphi=E$ or $\varphi=-E$&$\psi=E$ or $\psi=-E$\\\cline{3-5}
~&${\rm det}(\psi)=1$&1.2 &$\varphi =
                \begin{pmatrix}
                1 + mp^{2}q & mpq^{2}\\
                -mp^{3} & 1 - mp^{2}q
                \end{pmatrix}$&$\psi =
                \begin{pmatrix}
                1 + mpq^{2} & mq^{3}\\
                -mp^{2}q & 1 - mpq^{2}
                \end{pmatrix}$\\\cline{4-5}
 ~&~&&\multicolumn{2}{l|}{
$p, q, m \in \mathbb{Z}$ and
                ${\rm gcd}(p,q) = 1$}\\\cline{3-5}
~&~&1.3 &$\varphi=E$&$\psi \equiv_{3}
                \begin{pmatrix}
                1 & *\\
                0 & 1
                \end{pmatrix}$, ${\rm tr}(\psi) = -1$\\\cline{3-5}
~&~&1.4&$\varphi \equiv_{3} 
                \begin{pmatrix}
                1 & 0\\
                * & 1
                \end{pmatrix}$, ${\rm tr}(\varphi) = -1$&$\psi=E$\\\cline{3-5}
~&~&1.5 &$\varphi \equiv_{3} 
                \begin{pmatrix}
                0 & 2\\
                1 & 2
                \end{pmatrix}$ or
                $\varphi \equiv_{3} 
                \begin{pmatrix}
                2 & 1\\
                2 & 0
                \end{pmatrix}$ or
                $\varphi \equiv_{3} 
                E$, ${\rm tr}(\varphi)=-1$&$\psi=\varphi$\\\cline{3-5}
~&~&1.6 &$\varphi \equiv_{3} 
                \begin{pmatrix}
                0 & 1\\
                2 & 2
                \end{pmatrix}$  or $\varphi \equiv_{3} 
                \begin{pmatrix}
                2 & 2\\
                1 & 0
                \end{pmatrix}$ or $\varphi \equiv_{3} E$, ${\rm tr}(\varphi)=-1$&$\psi=\varphi^{-1}$\\
\hline
2.&${\rm det}(\varphi)=1$&2.1&$\varphi=E$&$\psi \equiv_{2} 
                \begin{pmatrix}
                1 & *\\
                0 & 1
                \end{pmatrix}, {\rm tr}(\psi) = 0$\\\cline{3-5}
~&${\rm det}(\psi)=-1$&2.2&$\varphi=-E$&$\psi \equiv_{2} E, {\rm tr}(\psi) = 0$\\
\hline
3.&${\rm det}(\varphi)=-1$&3.1&$\varphi \equiv_{2} 
            \begin{pmatrix}
            1 & 0\\
            * & 1
            \end{pmatrix}, {\rm tr}(\varphi) = 0$&$\psi=E$\\\cline{3-5}
~&${\rm det}(\psi)=1$&3.2&$\varphi \equiv_{2}E, {\rm tr}(\varphi) = 0$&$\psi=-E$\\
\hline
4.&${\rm det}(\varphi)=-1$&4.1&    $\varphi \equiv_{2} E$ or $\varphi \equiv_{2} 
                \begin{pmatrix}
                0 & 1\\
                1 & 0
                \end{pmatrix}$, ${\rm tr}(\varphi) = 0$&$\psi=\varphi$\\\cline{3-5}
~&${\rm det}(\psi)=-1$&4.2&  $\varphi \equiv_{2}E$, ${\rm tr}(\varphi) = 0$&$\psi=-\varphi$\\
\hline
\end{longtable}
}
\begin{center}
Table 1: Pairs of matrices $(\varphi,\psi)$ which lead to $\lambda$-homomorphic skew braces with $A_{\oplus}=\mathbb{Z}^2$.
\end{center}

Matrices on the left of formula (\ref{bardakovresultsyouknow}) are listed in Table~1 under case 1.2 for $p=q=1$, and matrices on the right of formula (\ref{bardakovresultsyouknow}) are listed in Table~1 under case 4.1. Hence the main result of the present paper generalizes the result of Bardakov, Neshchadim and Yadav \cite{BarNesYad}.

\section{Preliminaries}\label{prelim}
In this section, we recall the necessary definitions and formulate the statements which we are going to use throughout the paper.

Let $A$ be a skew brace.  For $a\in A$ denote by $\lambda_a$ the map $\lambda_a: x\mapsto \ominus a\oplus (a\odot x)$. Equality~(\ref{mainleft}) implies that $\lambda_a$ is an automorphism of $A_{\oplus}$. Moreover, the map $\lambda:A_{\odot}\to {\rm Aut}(A_{\oplus})$ given by $a\mapsto \lambda_a$ for $a\in A$ is a homomorphism. From the definition of $\lambda_a$ it follows that for all $a,b\in A$ the equality
$$
a\odot b=a\oplus \lambda_a(b)
$$
holds. Hence in order to define a skew brace it is sufficient to define the operation $\oplus$ and the maps $\lambda_a$ for all $a\in A$ appropriately. 

\subsection*{Skew braces and the holomorph}
Recall that the holomorph of a group $G$ is the semidirect product
$${\rm Hol}(G)=G\rtimes {\rm Aut}(G)=\{(g,f)~|~g\in G, f\in {\rm Aut}(G)\},$$ 
where ${\rm Aut}(G)$ acts on $G$ in a natural way. Denote by $\pi: {\rm Hol}(G)\to G$ the projection map given by the formula $\pi((g,f))=g$. A subgroup $H\leq {\rm Hol}(G)$ is called regular if for every element $g\in G$ there exists a unique element $(h,f)\in H$ such that $hf(g)=1$. It is known that if $H$ is a regular subgroup of ${\rm Hol}(G)$, then the map $\pi:H\to G$ is bijective. The following result  proved in \cite{Bac} establishes a one-to-one correspondence between skew braces with a given group $A_{\oplus}=G$ and regular subgroups of ${\rm Hol}(G)$. 
\begin{lemma}\label{th1}Let $G=(A,\oplus)$ be a group, and $A=(A,\oplus,\odot)$ be a skew brace. Then the group 
$$H_{\lambda}(A)=\{(a,\lambda_a)~|~a\in A\}$$
 is a regular subgroup of ${\rm Hol}(G)$. Moreover, the map $H_{\lambda}$ which sends $A$ to $H_{\lambda}(A)$ gives a bijection between skew braces $A=(A,\oplus,\odot)$ with  $A_{\oplus}=G$ and regular subgroups of ${\rm Hol}(G)$. 
\end{lemma}
Thus, in order to enumerate all skew braces with a given additive group $A_{\oplus}=G$ it suffices to enumerate all regular subgroups of ${\rm Hol}(G)$. If $H$ is a regular subgroup of ${\rm Hol}(G)$, and $(h,f)$ is an element of $H$, then $f=\lambda_h$, i.~e. $H$ gives the information about both  the operation $\oplus$ and the maps $\lambda_a$ for all $a\in G$.

\subsection*{$\lambda$-homomorphic skew braces}
We know that for every skew brace $A=(A,\oplus,\odot)$ the map $\lambda:a\mapsto \lambda_a$ is a homomorphism from $A_{\odot}$ to ${\rm Aut}(A_{\oplus})$. A skew brace is said to be $\lambda$-homomorphic, if this map is also a homomorphism from $A_{\oplus}$ to ${\rm Aut}(A_{\oplus})$.  If $A=(A,\oplus,\odot)$ is a $\lambda$-homomorphic skew braces, then for all $a,b\in A$ we have
\begin{align*}
\lambda_{a\oplus b}=\lambda_a\lambda_b,&&\lambda_{\ominus a}=\lambda_a^{-1}.
\end{align*}
Therefore, in order to definte $\lambda_a$ for all $a\in A$ it is enough to define $\lambda_a$ for all $a\in X$, where $X$ is the generating set of $A_{\oplus}$. The following statements are proved in \cite{BarNesYad}.
\begin{lemma}\label{th2} Let $A_{\oplus}=(A, \oplus)$ be a group, $X\subseteq A$ be the generating set of $A$, and $\lambda: A_{\oplus} \to {\rm Aut}(A_{\oplus})$ be a homomorphism with $\lambda(a)=\lambda_a$. Then the set 
 $$H_{\lambda} = \{(a,\lambda_a)~|~ a \in A\}$$ 
is a subgroup of ${\rm Hol}(A_{\oplus})$ if and only if for all $x, y \in X$ the element $\ominus x \oplus \lambda_{y}(x)$ belongs to the kernel of $ \lambda$. This subgroup is regular. \end{lemma}

The following result follows from Lemmas~\ref{th1} and Lemma \ref{th2}.
\begin{lemma}\label{th3} Let $A_{\oplus}=(A, \oplus)$ be a group, $X\subseteq A$ be the generating set of $A_{\oplus}$, and for each $x\in X$ the automorphism $\lambda_x\in{\rm Aut}(A_{\oplus})$ is defined. If the following two conditions hold
\begin{enumerate}
\item The map $x\mapsto \lambda_x$ for all $x\in X$ induces a homomorphism $\lambda:A_{\oplus}\to {\rm Aut}(A_{\oplus})$,
\item For all $x, y \in X$ the element $\ominus x \oplus \lambda_{y}(x)$ belongs to the kernel of $ \lambda$,
\end{enumerate}
then the algebraic system $(A,\oplus,\odot)$ with $a\odot b=a\oplus\lambda_a(b)$ for $a,b\in A$ is a $\lambda$-homomorphic skew brace.

Conversely, let $(A,\oplus,\odot)$ be a $\lambda$-homomorphic skew brace, such that the group $(A,\oplus)$ is generated by the set $X$. For $x\in X$ denote by $\lambda_x$ the automorphism which acts as $\lambda_x(a)=\ominus x\oplus (x\odot a)$. Then the two conditions specified above hold. 
\end{lemma}
\noindent \textbf{Proof.} If the two specified conditions hold, then from Lemma~\ref{th2} it follows that 
 $$H_{\lambda} = \{(a,\lambda_a)~|~ a \in A\}$$ 
is a regular subgroup of ${\rm Hol}(A_{\oplus})$. From Lemma~\ref{th1} it follows that this subgroup leads to the specified skew brace $A$. Since $\lambda:A_{\oplus}\to {\rm Aut}(A_{\oplus})$ is homomorphism, this brace is $\lambda$-homomorphic.

Conversely, if $(A,\oplus,\odot)$ is a $\lambda$-homomorphic skew brace, then the map $x\mapsto\lambda_x$ for all $x\in X$ induces a homomorphism $A_{\oplus}\to {\rm Aut}(A_{\oplus})$, and the first of two specified conditions hold. From Lemma~\ref{th1} it follows that 
$$H_{\lambda} = \{(a,\lambda_a)~|~ a \in A\}$$ 
is a subgroup of ${\rm Hol}(A_{\oplus})$, and from Lemma~\ref{th2} it follows that the second of two specified conditions hold.\hfill$\square$

Lemma~\ref{th3} implies that in order to classify all $\lambda$-homomorphic skew braces $(A,\oplus,\odot)$ with a given group $(A,\oplus)$ generated by a set $X$, it is enough to classify all sets $\{\lambda_x~|~x\in X\}$, which satisfy two conditions specified in Lemma~\ref{th3}.

\subsection*{Orders and centralizers of integer matrices}
Denote by $E$ the identity matrix in ${\rm GL}_2(\mathbb{Z})$. For $A\in {\rm GL}_2(\mathbb{Z})$ denote by $|A|$ the multiplicative order of $A$, i.~e. the minimal integer $n$ such that $A^n=E$. If such an integer doesn't exist, then assume that $|A|=\infty$. It is a well known fact that if  $A$ is a matrix from ${\rm GL}_2(\mathbb{Z})$, then $|A|\in \{1,2,3,4,6,\infty\}$. The following statement which characterizes matrices of finite order is proved, for example, in \cite{Baake-Roberts1997}.
\begin{lemma}\label{orders}Let $A$ be a matrix from ${\rm GL}_2(\mathbb{Z})$. Then the following statements hold.
\begin{enumerate}
    \item $|A| = 2$ if and only if $A=-E$ or ${\rm det}(A) = -1$ and ${\rm tr}(A) = 0$.
 \item $|A| = 3$ if and only if ${\rm det}(A) = 1$ and ${\rm tr}(A) = -1$.
 \item $|A| = 4$ if and only if ${\rm det}(A) = 1$ and ${\rm tr}(A) = 0$.
 \item $|A| = 6$ if and only if ${\rm det}(A) = 1$ and ${\rm tr}(A) = 1$.
\end{enumerate}
\end{lemma}
For a matrix $A$ from ${\rm GL}_2(\mathbb{Z})$ denote by $C(A)$ the centralizer of $A$ in ${\rm GL}_2(\mathbb{Z})$. The following statement is proved in \cite[Theorem 1]{Baake-Roberts1997}.
\begin{lemma}\label{centralizers}Let $A$ be a matrix from ${\rm GL}_2(\mathbb{Z})$. Then the following statements hold.
    \begin{enumerate}
        \item If $|A|=2$ or $|A|=4$, and $A\neq -E$, then $C(A) = \{ \pm E, \pm A\}$.
        \item If $|A|=3$ or $|A|=6$, then $C(A) = \{ \pm E, \pm A, \pm A^{-1}\}$.
        \item If $|A|= \infty$, then there exists  $B \in {\rm GL}_{2}(\mathbb{Z})$ such that $|B|=\infty$ and $C(A) = \{\pm B^m~|~m \in \mathbb{Z}\}$.
    \end{enumerate}
\end{lemma}

\section{$\lambda$-homomorphic braces on $\mathbb{Z}^2$}
Let 
$A_{\oplus}=(A,\oplus)=\mathbb{Z}^2=\langle x,y~|~[x,y]=1\rangle$
be the free abelian group with the generators $x,y$. Associate an automorphism $\varphi\in{\rm Aut}(A_{\oplus})$ to the element $x$, and an automorphism $\psi\in{\rm Aut}(A_{\oplus})$ to the element $y$. Since  ${\rm Aut}(A_{\oplus})={\rm Aut}(\mathbb{Z}^2)={\rm GL}_2(\mathbb{Z})$, we have
\begin{align*}
\varphi =
\begin{pmatrix}
\varphi_{11} & \varphi_{12}\\
\varphi_{21} & \varphi_{22}
\end{pmatrix},&&\psi =
\begin{pmatrix}
\psi_{11} & \psi_{12}\\
\psi_{21} & \psi_{22}
\end{pmatrix},
\end{align*}
where $\varphi_{i,j}, \psi_{i,j}\in \mathbb{Z}$, and ${\rm det}(\varphi), {\rm det}(\psi)\in \{1,-1\}$. The following statement says when the automorphisms $\varphi,\psi$ satisfy two conditions specified in Lemma~\ref{th3}. 
\begin{lemma}\label{conclemma} Let 
$A_{\oplus}=(A,\oplus)=\mathbb{Z}^2$ be the free abelian group with the generators $x,y$. Associate an automorphism $\varphi\in{\rm Aut}(A_{\oplus})$ to $x$, and an automorphism $\psi\in{\rm Aut}(A_{\oplus})$ to $y$. Then the two conditions specified in Lemma~\ref{th3} hold if an only if $\varphi\psi=\psi\varphi$ and the following equalities hold
\begin{align}\label{main equalities}
\varphi^{\varphi_{11}-1}\psi^{\varphi_{21}} = E,&&
\varphi^{\varphi_{12}}\psi^{\varphi_{22}-1} = E,&&
\varphi^{\psi_{11}-1}\psi^{\psi_{21}} = E,&&
\varphi^{\psi_{12}}\psi^{\psi_{22}-1} = E.
\end{align}
\end{lemma}
\noindent \textbf{Proof.} The first condition specified in Lemma~\ref{th3} says that the map $x\mapsto \lambda_x$ for all $x\in X$ induces a homomorphism $\lambda:A_{\oplus}\to {\rm Aut}(A_{\oplus})$. It means that the maps $x\mapsto \varphi$, $y\mapsto \psi$ induce a homomorphism $\mathbb{Z}^2\to {\rm GL}_2(\mathbb{Z})$. It happens if and only if $\varphi\psi=\psi\varphi$. 

The second condition specified in Lemma~\ref{th3} says that for all $x, y \in X$ the element $\ominus x \oplus \lambda_{y}(x)$ belongs to the kernel of $ \lambda$. This condition can be written in the following form.
\begin{align*}
\ominus x\oplus\varphi(x)  \in {\rm Ker}(\lambda) &\Leftrightarrow \varphi^{\varphi_{11}-1}\psi^{\varphi_{21}} = E\\
\ominus y\oplus\varphi(y)  \in {\rm Ker}(\lambda) &\Leftrightarrow \varphi^{\varphi_{12}}\psi^{\varphi_{22}-1} = E\\
\ominus x\oplus\psi(x)  \in {\rm Ker}(\lambda) &\Leftrightarrow\varphi^{\psi_{11}-1}\psi^{\psi_{21}} = E\\
\ominus y\oplus\psi(y)  \in {\rm Ker}(\lambda) &\Leftrightarrow \varphi^{\psi_{12}}\psi^{\psi_{22}-1} = E
\end{align*}
Lemma is proved.\hfill$\square$

Thus in order to classify all $\lambda$-homomorphic skew braces $(A,\oplus,\odot)$ with $A_{\oplus}=\mathbb{Z}^2$ it is enough to find all automorphisms $\varphi,\psi\in{\rm GL}_2(\mathbb{Z})$ such that $\varphi\psi=\psi\varphi$ and equalities (\ref{main equalities}) hold. The following theorem describes such automorphisms. 
\begin{theorem}\label{ttthhhtt}Let the automorphisms $\varphi,\psi\in{\rm GL}_2(\mathbb{Z})$ be such that $\varphi\psi=\psi\varphi$ and equalities (\ref{main equalities}) hold.  Then  $(\varphi,\psi)$ satisfy one of the conditions listed in Table 1.
\end{theorem} 

\noindent \textbf{Proof.} It is clear that the matrices $\varphi,\psi\in\{E,-E\}$ satisfy the formulation of the theorem. Hence, in the proof we assume that at least one of matrices $\varphi,\psi$ is not equal to $\pm E$. Moreover, we will consider in details only the case when $\varphi\neq\pm E$. The case when $\varphi=\pm E$ and $\psi \neq \pm E$ is similar, and we will mention this case without details. From Section~\ref{prelim} we know that the order of $\varphi\neq \pm E$ belongs to $\{2,3,4,6,\infty\}$. Depending on $|\varphi|$, consider several cases.

\textbf{Case 1: $|\varphi|=2$ and $\varphi\neq -E$.} From the condition $\varphi\psi=\psi\varphi$ it follows that $\psi$ belongs to the centralizer $C(\varphi)$. From Lemma~\ref{centralizers} it follows that $C(\varphi)=\{\pm E, \pm \varphi\}$. Depending on $\psi$ consider several cases.

If $\psi=E$, then conditions~(\ref{main equalities}) can be rewritten in the form $\varphi^{\varphi_{11}-1}=E$, $\varphi^{\varphi_{12}}=E$. These conditions mean that $\varphi_{11}\equiv_21$, $\varphi_{12}\equiv_20$ or
$$\varphi=\begin{pmatrix}1+2m&2p\\
q&n
\end{pmatrix}$$
for appropriate $p,q,m,n\in\mathbb{Z}$. Since $|\varphi|=2$ and $\varphi\neq -E$, from Lemma~\ref{orders} we have ${\rm tr}(\varphi)=0$ and ${\rm det}(\varphi)=-1$. It means that $n=-1-2m$ and 
$$-1={\rm det}(\varphi)={\rm det}\begin{pmatrix}1+2m&2p\\
q&-1-2m
\end{pmatrix}=-(1+2m)^2 - 2pq.$$
This equation has two solutions
\begin{align*}
m = \dfrac{-1\pm\sqrt{1-2pq}}{2}.
\end{align*}
Since $1$ is a quadratic deduction modulo $2$, parameter $m$ in the last equation can be integer-valued for appropriate integers $pq$. Hence, $\varphi$ has the following form
$$\varphi=\begin{pmatrix}
\pm\sqrt{1-2pq} & 2p\\
q & \mp\sqrt{1-2pq}
\end{pmatrix}$$
for arbitrary $p,q\in \mathbb{Z}$ which provide integer values in the entries of $\varphi$. This is exactly case 3.1 from Table 1. The case $\psi=E$ is finished. The symmetric case ($\varphi=E$, $|\psi|=2$, $\psi\neq \pm E$) is case 2.1 from Table 1.

If $\psi=-E$, then conditions~(\ref{main equalities}) can be rewritten in the form $\varphi^{\varphi_{11}-1}(-E)^{\varphi_{21}} = E$, $
\varphi^{\varphi_{12}}(-E)^{\varphi_{22}-1} = E$. Since $|\varphi|=2$ and $\varphi\neq -E$, there is no power of $\varphi$ which is equal to $-E$. Therefore the last system of equalities can be rewritten in the form
\begin{align*}
\varphi_{11}\equiv_21,&&\varphi_{21}\equiv_20,&&\varphi_{12}\equiv_20,&&\varphi_{22}\equiv_21.
\end{align*}
It means that
$$\varphi=\begin{pmatrix}1+2m&2p\\
2q&1+2n
\end{pmatrix}$$
for appropriate $p,q,m,n\in\mathbb{Z}$. Since $|\varphi|=2$ and $\varphi\neq -E$, from Lemma~\ref{orders} we have ${\rm tr}(\varphi)=0$ and ${\rm det}(\varphi)=-1$. It means that $n=-1-m$ and 
$$-1={\rm det}(\varphi)={\rm det}\begin{pmatrix}1+2m&2p\\
2q&-1-2m
\end{pmatrix}=-(1+2m)^2 - 4pq.$$
This equation has two solutions
\begin{align*}
m = \dfrac{-1\pm\sqrt{1-4pq}}{2}.
\end{align*}
Since $1$ is a quadratic deduction modulo $4$, parameter $m$ in the last equation can be integer-valued for appropriate integers $pq$. Hence, $\varphi$ has the following form
$$\varphi=\begin{pmatrix}
\pm\sqrt{1-4pq} & 2p\\
2q & \mp\sqrt{1-4pq}
\end{pmatrix}$$
for arbitrary $p,q\in \mathbb{Z}$ which provide integer values in the entries of $\varphi$. This is exactly case 3.2 from Table 1. The case $\psi=-E$ is finished. The symmetric case ($\varphi=-E$, $|\psi|=2$, $\psi\neq \pm E$) is case 2.2 from Table 1.

If $\psi=\varphi$, then conditions~(\ref{main equalities}) can be rewritten in the form $\varphi^{\varphi_{11}+\varphi_{21}-1} = E$, $\varphi^{\varphi_{12}+\varphi_{22}-1}= E$. These conditions mean that $\varphi_{11} + \varphi_{21} \equiv_2 1$, $\varphi_{12} + \varphi_{22} \equiv_2 1$. Therefore $\varphi$ can be written as
$$\varphi=\begin{pmatrix}
\varphi_{11} & \varphi_{12}\\
1+2m-\varphi_{11} & 1+2n-\varphi_{12}   
\end{pmatrix}$$
for appropriate $m,n\in\mathbb{Z}$. Since $|\varphi|=2$ and $\varphi\neq -E$, from Lemma~\ref{orders} we have ${\rm tr}(\varphi)=0$ and ${\rm det}(\varphi)=-1$. From equality ${\rm tr}(\varphi)=0$ we have $\varphi_{12} = 1 + 2n + \varphi_{11}$. Therefore the equality ${\rm det}(\varphi)=-1$ can be written as
$$-1={\rm det}(\varphi)={\rm det}\begin{pmatrix}
\varphi_{11} & 1 + 2n + \varphi_{11} \\
1+2m-\varphi_{11} & -\varphi_{11}   
\end{pmatrix}=2\varphi_{11}(n-m)-4mn-2(m+n)-1.$$
This equality means that 
\begin{equation}\label{one of the first equations}(m-n)\varphi_{11} + 2mn+m+n = 0.
\end{equation}
If $n=m$, then this equality can be rewritten as $m^2+m=0$. Hence, we have either $m=n=0$ or $m=n=-1$, and if we denote by $\varphi_{11}=p$, then we have one of the following cases
\begin{align*}
\varphi = \begin{pmatrix}
p & 1 + p\\
1-p & -p   
\end{pmatrix},&&
\varphi = \begin{pmatrix}
p & -1 + p\\
-1 - p & -p   
\end{pmatrix}.
\end{align*}
These cases correspond to case 4.1 from Table 1. 

If $n\neq m$, then equation (\ref{one of the first equations}) implies that
$\varphi_{11}=\dfrac{m+n+2mn}{n-m}$
and
$$\varphi = \begin{pmatrix}
\dfrac{m+n+2mn}{n-m} & \dfrac{2n(1+n)}{n-m}\\
\dfrac{-2m(1+m)}{n-m} & \dfrac{-m-n-2mn}{n-m} 
\end{pmatrix}.$$
For appropriate $p,q\in \mathbb{Z}$ this matrix can have integer entries. In this case this matrix corresponds  to case 4.1 from Table 1.

If $\psi=-\varphi$, then conditions~(\ref{main equalities}) can be rewritten in the form
\begin{align*}
\varphi^{\varphi_{11}+\varphi_{21}-1}(-E)^{\varphi_{21}} &= E,\\
\varphi^{\varphi_{12}+\varphi_{22}-1}(-E)^{\varphi_{22}-1} &= E,\\
\varphi^{-\varphi_{11}-\varphi_{21}-1}(-E)^{-\varphi_{21}} &= E,\\
\varphi^{-\varphi_{12}-\varphi_{22}-1}(-E)^{-\varphi_{22}-1} &= E.
\end{align*}
Since $|\varphi|=2$ and $\varphi\neq -E$, there is no power of $\varphi$ which is equal to $-E$. Therefore the last system of equalities can be rewritten in the form
\begin{align*}
\varphi_{11} + \varphi_{21} \equiv_2 1,&& \varphi_{12} + \varphi_{22} \equiv_2 1,&&\varphi_{21}\equiv_20,&&\varphi_{22}\equiv_21.
\end{align*}
This means that
$$
\varphi=\begin{pmatrix}1+2n&2m\\
2p&1+2q
\end{pmatrix}.
$$
Since $|\varphi|=2$ and $\varphi\neq -E$, from Lemma~\ref{orders} we have ${\rm tr}(\varphi)=0$ and ${\rm det}(\varphi)=-1$. From equality ${\rm tr}(\varphi)=0$ we have $q = -1 -n$. Therefore the equality ${\rm det}(\varphi)=-1$ can be written as
$$-1={\rm det}(\varphi)={\rm det}\begin{pmatrix}
1+2n & 2m \\
2p & -1-2n   
\end{pmatrix}=-(1+2n)^2-4mp.$$
This equation has two solutions
\begin{align*}
n = \dfrac{-1\pm\sqrt{1-4mp}}{2}.
\end{align*}
Since $1$ is a quadratic deduction modulo $4$, parameter $n$ in the last equation can be integer-valued for appropriate integers $m, p$. Hence, $\varphi$ has the following form
$$\varphi=\begin{pmatrix}
\pm\sqrt{1-4mp} & 2m\\
2p & \mp\sqrt{1-4mp}
\end{pmatrix}$$
for arbitrary $m,p\in \mathbb{Z}$ which provide integer values in the entries of $\varphi$. This is exactly case 4.2 from Table 1. The case $\psi=-\varphi$ is finished, and the case $|\varphi|=2$, $\varphi\neq -E$ is finished.

\textbf{Case 2: $|\varphi|=3$.} From the condition $\varphi\psi=\psi\varphi$ it follows that $\psi$ belongs to the centralizer $C(\varphi)$. From Lemma~\ref{centralizers} it follows that $C(\varphi)=\{\pm E, \pm \varphi, \pm \varphi^{-1}\}$. Depending on $\psi$ consider several cases.

If $\psi=E$, then conditions~(\ref{main equalities}) can be rewritten in the form $\varphi^{\varphi_{11}-1}=E$, $\varphi^{\varphi_{12}}=E$.
These conditions mean that $\varphi_{11}\equiv_31$, $\varphi_{12}\equiv_30$ or
$$\varphi=\begin{pmatrix}1+3m&3p\\
q&n
\end{pmatrix}$$
for appropriate $p,q,m,n\in\mathbb{Z}$. Since $|\varphi|=3$, from Lemma~\ref{orders} we have ${\rm tr}(\varphi)=-1$ and ${\rm det}(\varphi)=1$. It means that $n=-2-3m$ and 
$$1={\rm det}(\varphi)={\rm det}\begin{pmatrix}1+3m&3p\\
q&-2-3m
\end{pmatrix}=-(2+3m)(1+3m) - 3pq.$$
The last equality can be rewritten as the equation $3m^2+3m+pq+1=0$ over variable $m$. This equation has two solutions
\begin{align*}
m = \dfrac{-3\pm\sqrt{-3-12pq}}{6}.
\end{align*}
Since $9$ is a quadratic deduction modulo $12$, parameter $m$ in the last equation can be integer-valued for appropriate integers $pq$. Hence, $\varphi$ has the following form
$$\varphi=\begin{pmatrix}
\dfrac{-1\pm\sqrt{-3-12pq}}{2} & 3p\\
q & \dfrac{-1\mp\sqrt{-3-12pq}}{2}
\end{pmatrix}$$
for arbitrary $p,q\in \mathbb{Z}$ which provide integer values in the entries of $\varphi$. This is exactly case 1.4 from Table 1. The case $\psi=E$ is finished. The symmetric case ($\varphi=E$, $|\psi|=3$) is case 1.3 from Table 1.

If $\psi=-E$, then the third equality from conditions~(\ref{main equalities}) can be rewritten in the form $\varphi^{-2} = E$. Since $|\varphi|=3$, this equality cannot hold, and the case $\psi=-E$ is finished.

If $\psi=\varphi$, then conditions~(\ref{main equalities}) can be rewritten in the form $\varphi^{\varphi_{11}+\varphi_{21}-1} = E$, $
\varphi^{\varphi_{12}+\varphi_{22}-1}= E$. These conditions mean that $\varphi_{11} + \varphi_{21} \equiv_3 1$, $\varphi_{12} + \varphi_{22} \equiv_3 1$. Therefore $\varphi$ can be written as
$$\varphi=\begin{pmatrix}
\varphi_{11} & \varphi_{12}\\
1+3m-\varphi_{11} & 1+3n-\varphi_{12}   
\end{pmatrix}$$
for appropriate $m,n\in\mathbb{Z}$. Since $|\varphi|=3$, from Lemma~\ref{orders} we have ${\rm tr}(\varphi)=-1$ and ${\rm det}(\varphi)=1$. From equality ${\rm tr}(\varphi)=-1$ we have $\varphi_{12} = 2 + 3n + \varphi_{11}$. Therefore the equality ${\rm det}(\varphi)=1$ can be written as
$$1={\rm det}(\varphi)={\rm det}\begin{pmatrix}
\varphi_{11} & 2 + 3n + \varphi_{11} \\
1+3m-\varphi_{11} & -1-\varphi_{11}   
\end{pmatrix}=3\varphi_{11}(n-m)-9mn-6m-3n-2.$$
This equality means that $(m-n)\varphi_{11} + n + 2m + 3mn + 1 = 0$. If $n=m$, then this equality can be rewritten as $3m^2+3m+1=0$. This quadratic equation doesn't have integer solutions. If $n\neq m$, then 
\begin{equation}\label{f11}
\varphi_{11} = \dfrac{1+n+2m+3mn}{n-m},
\end{equation}
and for appropriate $m,n\in \mathbb{Z}$ the number $\varphi_{11}$ can belong to $\mathbb{Z}$. Hence, $\varphi$ is of the form 
$$\begin{pmatrix}
\varphi_{11} & 2 + 3n + \varphi_{11} \\
1+3m-\varphi_{11} & -1-\varphi_{11}   
\end{pmatrix},$$
where $\varphi_{11}$ is given by formula~(\ref{f11}), and $n\neq m$ are such integers that $\varphi_{11}$ is also an integer. This is exactly case 1.5 from Table 1. The case $\psi=\varphi$ is finished. 

If $\psi=-\varphi$, then conditions~(\ref{main equalities}) can be rewritten in the form
\begin{align}
\label{eq1ord3}\varphi^{\varphi_{11}+\varphi_{21}-1}(-E)^{\varphi_{21}} &= E,\\
\varphi^{\varphi_{12}+\varphi_{22}-1}(-E)^{\varphi_{22}-1} &= E,\\
\label{eq3ord3}\varphi^{-\varphi_{11}-\varphi_{21}-1}(-E)^{-\varphi_{21}} &= E,\\
\varphi^{-\varphi_{12}-\varphi_{22}-1}(-E)^{-\varphi_{22}-1} &= E.
\end{align}
Multiplying equalities (\ref{eq1ord3}) and (\ref{eq3ord3}) we conculde that $\varphi^{-2}=E$, what contradicts the fact that $|\varphi|=3$. So, the situation $\psi=-\varphi$ is impossible, and the case $\psi=-\varphi$ is finished.

If $\psi=\varphi^{-1}$, then 
$$\psi=\begin{pmatrix}
\varphi_{22} & -\varphi_{12} \\
-\varphi_{21} & \varphi_{11} 
\end{pmatrix}$$
and  conditions~(\ref{main equalities}) can be rewritten in the form
\begin{align*}
\varphi^{\varphi_{11}-\varphi_{21}-1} = E,&&
\varphi^{\varphi_{12}-\varphi_{22}+1} = E,&&
\varphi^{\varphi_{22}+\varphi_{21}-1} = E,&&
\varphi^{-\varphi_{12}-\varphi_{11}+1} = E.
\end{align*}
These equalities mean that,
\begin{equation}\label{newsystemequiv}\begin{cases}\varphi_{11}-\varphi_{21}\equiv_31,\\
\varphi_{12}-\varphi_{22}\equiv_3-1,\\
\varphi_{22}+\varphi_{21}\equiv_31,\\
-\varphi_{12}-\varphi_{11}\equiv_3-1.
\end{cases}
\end{equation}
From Lemma~\ref{orders} we have that ${\rm tr}(\varphi)=-1$. Using this equality, system (\ref{newsystemequiv}) can be written as
$$\begin{cases}
\varphi_{11} - \varphi_{21} \equiv_{3} 1,\\ 
\varphi_{11} + \varphi_{12} \equiv_{3} 1.
\end{cases}
$$
It means that
\begin{equation}\label{phiinthiscase}
\varphi =\begin{pmatrix}
\varphi_{11} & 1+3n-\varphi_{11}\\
-1-3m+\varphi_{11} & -1-\varphi_{11}  
\end{pmatrix}
\end{equation}
for appropriate integers $m,n$. From Lemma~\ref{orders} we have ${\rm det}(\varphi)=1$, what means that 
\begin{equation}\label{onemorem1111}
(1+m+n)\varphi_{11} - (m+n) -3mn = 0.
\end{equation}

If $1 + n + m = 0$, then $m = -n -1$, and equality~(\ref{onemorem1111}) can be rewritten as $3n^2 + 3n + 1 = 0$. This quadratic equation doesn't have integer solutions. If $1 + n + m \neq 0$, then from equality (\ref{onemorem1111}) we have
\begin{equation}\label{ord3var5}
\varphi_{11} = \dfrac{3mn+m+n}{1+n+m},
\end{equation}
and $\varphi_{11}$ can be an integer for appropriate $m,n\in \mathbb{Z}$. Hence, $\varphi$ is of the form (\ref{phiinthiscase}), where $\varphi_{11}$ is given by formula~(\ref{ord3var5}), and $n,m$ are such integers that $\varphi_{11}$ is also an integer. This is exactly case 1.6 from Table 1. The case $\psi=\varphi^{-1}$ is finished. 

If $\psi=-\varphi^{-1}$, then 
$$\psi=\begin{pmatrix}
-\varphi_{22} & \varphi_{12} \\
\varphi_{21} & -\varphi_{11} 
\end{pmatrix}$$
and  conditions~(\ref{main equalities}) can be rewritten in the form
\begin{align}
\label{qqqeee1}\varphi^{\varphi_{11}-\varphi_{21}-1}(-E)^{\varphi_{21}} &= E,\\
\varphi^{\varphi_{12}-\varphi_{22}+1}(-E)^{\varphi_{22}-1} &= E,\\
\label{qqqeee2}\varphi^{-\varphi_{22}-\varphi_{21}-1}(-E)^{\varphi_{21}} &= E,\\
\varphi^{\varphi_{12}+\varphi_{11}+1}(-E)^{-\varphi_{11}-1} &= E.
\end{align}
From Lemma~\ref{orders} we have that ${\rm tr}(\varphi)=-1$, hence, $\varphi_{22}=-1-\varphi_{11}$. Using this equality, equalities (\ref{qqqeee1}), (\ref{qqqeee2}) can be rewritten in the form
\begin{align*}
\varphi^{\varphi_{11}-\varphi_{21}-1}(-E)^{\varphi_{21}} = E,&&\varphi^{\varphi_{11}+1-\varphi_{21}-1}(-E)^{\varphi_{21}} = E.
\end{align*}
Multiplying the second of these equalities by the inverse to the first of these equalities we have the equality $\varphi=E$, which contradicts the equality $|\varphi|=3$. So, the situation $\psi=-\varphi^{-1}$ is impossible, and the case $\psi=-\varphi^2$ is finished. The case $|\varphi|=3$ is finished.

\textbf{Case 3: $|\varphi|=4$.} From the condition $\varphi\psi=\psi\varphi$ it follows that $\psi$ belongs to the centralizer $C(\varphi)$. From Lemma~\ref{centralizers} it follows that $C(\varphi)=\{\pm E, \pm \varphi\}$. Depending on $\psi$ consider several cases.

If $\psi=E$, then conditions~(\ref{main equalities}) can be rewritten in the form $\varphi^{\varphi_{11}-1}=E$, $\varphi^{\varphi_{12}}=E$. These conditions mean that $\varphi_{11}\equiv_41$, $\varphi_{12}\equiv_40$ or
$$\varphi=\begin{pmatrix}1+4m&4p\\
q&n
\end{pmatrix}$$
for appropriate $p,q,m,n\in\mathbb{Z}$. Since $|\varphi|=4$, from Lemma~\ref{orders} we have ${\rm tr}(\varphi)=0$ and ${\rm det}(\varphi)=1$. It means that $n=-1-4m$ and 
$$1={\rm det}(\varphi)={\rm det}\begin{pmatrix}1+4m&4p\\
q&-1-4m
\end{pmatrix}=-(1+4m)^2 - 4pq.$$
The last equality can be rewritten as the equation $8m^2+4m+1+2pq=0$ over variable $m$. This equation has two solutions
\begin{align*}
m = \dfrac{-1\pm\sqrt{-1-4pq}}{4}.
\end{align*}
Since $3$ is not a quadratic deduction modulo $4$, for any integers $p,q$ parameter $m$ in the last equation cannot be integer. 

If $\psi=-E$, then the third equality from conditions~(\ref{main equalities}) can be rewritten in the form $\varphi^{-2} = E$. Since $|\varphi|=4$, this equality cannot hold, and the case $\psi=-E$ is finished.

If $\psi=\varphi$, then conditions~(\ref{main equalities}) can be rewritten in the form $\varphi^{\varphi_{11}+\varphi_{21}-1} = E$, $\varphi^{\varphi_{12}+\varphi_{22}-1}= E$. These conditions mean that $\varphi_{11} + \varphi_{21} \equiv_4 1$, $\varphi_{12} + \varphi_{22} \equiv_4 1$. Therefore $\varphi$ can be written as
$$\varphi=\begin{pmatrix}
\varphi_{11} & \varphi_{12}\\
1+4m-\varphi_{11} & 1+4n-\varphi_{12}   
\end{pmatrix}$$
for appropriate $m,n\in\mathbb{Z}$. Since $|\varphi|=4$, from Lemma~\ref{orders} we have ${\rm tr}(\varphi)=0$ and ${\rm det}(\varphi)=1$. From equality ${\rm tr}(\varphi)=0$ we have $\varphi_{12} = 1 + 4n + \varphi_{11}$. Therefore the equality ${\rm det}(\varphi)=1$ can be written as
$$1={\rm det}(\varphi)={\rm det}\begin{pmatrix}
\varphi_{11} & 1 + 4n + \varphi_{11} \\
1+4m-\varphi_{11} & -\varphi_{11}   
\end{pmatrix}=4\varphi_{11}(n-m)-16mn-4(m+n)-1.$$
This equality means that  $2(n-m)\varphi_{11} -8mn-2(m+n)-1 = 0$. If $n=m$, then this equality can be rewritten as $8m^2+4m+1=0$. This quadratic equation doesn't have integer solutions. If $n\neq m$, then 
$$
\varphi_{11} = \dfrac{2(m+n+4mn)+1}{2(n-m)}
$$
is not an integer, and the case $\psi=\varphi$ is finished. 

If $\psi=-\varphi$, then repeating steps when $|\varphi|=3$ and $\psi=-\varphi$ we conclude that $\varphi^{-2}=E$, what contradicts the fact that $|\varphi|=4$. So, the situation $\psi=-\varphi$ is impossible, the case $\psi=-\varphi$ is finished, and the case $|\varphi|=4$ is finished.

\textbf{Case 4: $|\varphi|=6$.} From the condition $\varphi\psi=\psi\varphi$ it follows that $\psi$ belongs to the centralizer $C(\varphi)$. From Lemma~\ref{centralizers} it follows that $C(\varphi)=\{\pm E, \pm \varphi, \pm \varphi^{-1}\}$. Depending on $\psi$ consider several cases.

If $\psi=E$, then conditions~(\ref{main equalities}) can be rewritten in the form $\varphi^{\varphi_{11}-1}=E$, $\varphi^{\varphi_{12}}=E$. These conditions mean that $\varphi_{11}\equiv_61$, $\varphi_{12}\equiv_60$ or
$$\varphi=\begin{pmatrix}1+6m&6p\\
q&n
\end{pmatrix}$$
for appropriate $p,q,m,n\in\mathbb{Z}$. Since $|\varphi|=6$, from Lemma~\ref{orders} we have ${\rm tr}(\varphi)=1$ and ${\rm det}(\varphi)=1$. It means that $n=-6m$ and 
$$1={\rm det}(\varphi)={\rm det}\begin{pmatrix}1+6m&6p\\
q&-6m
\end{pmatrix}=-36m^2-6m-6pq.$$
The last equality can be rewritten as the equation $36m^2+6m+1+6pq=0$ over variable $m$. This equation has two solutions
\begin{align*}
m = \dfrac{-1\pm\sqrt{-3-24pq}}{12}.
\end{align*}
Since $3$ is not a quadratic deduction modulo $24$, for any integers $p,q$ parameter $m$ in the last equation cannot be integer. 

If $\psi=-E$, then the third equality from conditions~(\ref{main equalities}) can be rewritten in the form $\varphi^{-2} = E$. Since $|\varphi|=6$, this equality cannot hold, and the case $\psi=-E$ is finished.

If $\psi=\varphi$, then conditions~(\ref{main equalities}) can be rewritten in the form $\varphi^{\varphi_{11}+\varphi_{21}-1} = E$, $\varphi^{\varphi_{12}+\varphi_{22}-1}= E$. These conditions mean that $\varphi_{11} + \varphi_{21} \equiv_6 1$, $\varphi_{12} + \varphi_{22} \equiv_6 1$. Therefore $\varphi$ can be written as
$$\varphi=\begin{pmatrix}
\varphi_{11} & \varphi_{12}\\
1+6m-\varphi_{11} & 1+6n-\varphi_{12}   
\end{pmatrix}$$
for appropriate $m,n\in\mathbb{Z}$. Since $|\varphi|=6$, from Lemma~\ref{orders} we have ${\rm tr}(\varphi)=1$ and ${\rm det}(\varphi)=1$. From equality ${\rm tr}(\varphi)=1$ we have $\varphi_{12} = 6n + \varphi_{11}$. Therefore the equality ${\rm det}(\varphi)=1$ can be written as
$$1={\rm det}(\varphi)={\rm det}\begin{pmatrix}
\varphi_{11} & 6n + \varphi_{11} \\
1+6m-\varphi_{11} &1 -\varphi_{11}   
\end{pmatrix}=6(n-m)\varphi_{11}-6n-36mn.$$
If $n=m$, then this equality can be rewritten as $-6m(1+6m)=1$. This equation doesn't have integer solutions. If $n\neq m$, then 
$$
\varphi_{11} = \dfrac{6n(1+6m)+1}{6(n-m)}
$$
is not an integer, and the case $\psi=\varphi$ is finished. 

If $\psi=-\varphi$, then repeating steps when $|\varphi|=3$ and $\psi=-\varphi$ we conclude that $\varphi^{-2}=E$, what contradicts the fact that $|\varphi|=6$. So, the case $\psi=-\varphi$ is finished.

If $\psi=\varphi^{-1}$, then 
$$\psi=\begin{pmatrix}
\varphi_{22} & -\varphi_{12} \\
-\varphi_{21} & \varphi_{11} 
\end{pmatrix}$$
and  conditions~(\ref{main equalities}) can be rewritten in the form
\begin{align*}
\varphi^{\varphi_{11}-\varphi_{21}-1} = E,&&
\varphi^{\varphi_{12}-\varphi_{22}+1} = E,&&
\varphi^{\varphi_{22}+\varphi_{21}-1} = E,&&
\varphi^{-\varphi_{12}-\varphi_{11}+1} = E.
\end{align*}
These equalities mean that,
\begin{equation}\label{newsystemequiv2221}
\begin{cases}\varphi_{11}-\varphi_{21}\equiv_61,\\
\varphi_{12}-\varphi_{22}\equiv_6-1,\\
\varphi_{22}+\varphi_{21}\equiv_61,\\
-\varphi_{12}-\varphi_{11}\equiv_6-1.
\end{cases}
\end{equation}
From Lemma~\ref{orders} we have that ${\rm tr}(\varphi)=1$. Using this equality, system (\ref{newsystemequiv2221}) can be written as
$$\begin{cases}\varphi_{11}-\varphi_{21}\equiv_61,\\
\varphi_{12}+\varphi_{11}\equiv_60,\\
\varphi_{11}-\varphi_{21}\equiv_60,\\
+\varphi_{12}+\varphi_{11}\equiv_61.
\end{cases}
$$
This system doesn't have integer solutions.

If $\psi=-\varphi^{-1}$, then 
$$\psi=\begin{pmatrix}
-\varphi_{22} & \varphi_{12} \\
\varphi_{21} & -\varphi_{11} 
\end{pmatrix}$$
and  conditions~(\ref{main equalities}) can be rewritten in the form
\begin{align}
\label{qqqeee1444}\varphi^{\varphi_{11}-\varphi_{21}-1}(-E)^{\varphi_{21}} &= E,\\
\varphi^{\varphi_{12}-\varphi_{22}+1}(-E)^{\varphi_{22}-1} &= E,\\
\varphi^{-\varphi_{22}-\varphi_{21}-1}(-E)^{\varphi_{21}} &= E,\\
\label{qqqeee2444}\varphi^{\varphi_{12}+\varphi_{11}+1}(-E)^{-\varphi_{11}-1} &= E.
\end{align}

From Lemma~\ref{orders} we have that ${\rm tr}(\varphi)=1$, hence, $\varphi_{22}=1-\varphi_{11}$. Using this equality, equalities (\ref{qqqeee1444})- (\ref{qqqeee2444}) can be rewritten in the form
\begin{align*}
\varphi^{\varphi_{11}-\varphi_{21}-1}(-E)^{\varphi_{21}} &= E,\\
\varphi^{\varphi_{12}+\varphi_{11}}(-E)^{-\varphi_{11}} &= E,\\
\varphi^{\varphi_{11}-\varphi_{21}-2}(-E)^{\varphi_{21}} &= E,\\
\varphi^{\varphi_{12}+\varphi_{11}+1}(-E)^{-\varphi_{11}-1} &= E.
\end{align*}
Multiplying the first of these equalities by the inverse to the third of these equalities we have the equality $\varphi=E$, which contradicts the equality $|\varphi|=6$. So, the situation $\psi=-\varphi^{-1}$ is impossible, the case $\psi=-\varphi^{-1}$ is finished, and the case $|\varphi|=6$ is finished.

\textbf{Case 5: $|\varphi|=\infty$.}  From the condition $\varphi\psi=\psi\varphi$ it follows that $\psi$ belongs to the centralizer $C(\varphi)$. From Lemma~\ref{centralizers} it follows that $C(\varphi)=\{\pm \theta^m~|~m \in \mathbb{Z}\}$, where $\theta$ is an appropriate automorphism with $|\theta|=\infty$. Hence, $\varphi = \varepsilon_1\theta^p$, $\psi = \varepsilon_2\theta^q$, where $p,q\in \mathbb{Z}$ are appropriate integers, and $\varepsilon_{1},\varepsilon_2 = \pm 1$. 

If $q=0$, then $\psi=\varepsilon_2E$, and $\psi_{12}=\psi_{21}=0$. Therefore, equalities (\ref{main equalities}) can be rewritten in the following form
\begin{align}\label{we rewrite for infinite order}
\varphi^{\varphi_{11}-1}\varepsilon_2^{\varphi_{21}} = E,&&
\varphi^{\varphi_{12}}\varepsilon_2^{\varphi_{22}-1} = E,&&
\varphi^{\psi_{11}-1} = E,&&
\varepsilon_2^{\psi_{22}-1} = E.
\end{align}
Since $|\varphi|=\infty$, from the third of these equalities it follows that $\psi_{11}=1$. Hence $\psi=E$, and $\psi_{22}=1$.  Equalities (\ref{we rewrite for infinite order}) in this situation can be written as $\varphi^{\varphi_{11}-1} = E$, $\varphi^{\varphi_{12}} = E$.  Therefore $\varphi_{11}=1$, $\varphi_{12}=0$ and 
$$\varphi = 
\begin{pmatrix}
1  & 0\\
m  & 1
\end{pmatrix}.$$
This is the case 1.2 from Table 1 for $q=0, p=-1$.

If $q\neq 0$, then denote by $m=\gcd(p,q)$, $\theta_1=\theta^m$, $p_1=p/m$, $q_1=q/m$. In this notations we have $\varphi=\varepsilon_1\theta_1^{p_1}$, $\psi_2=\varepsilon_2\theta_1^{q_1}$, where $\gcd(p_1,q_1)=1$. Changing letters $\theta_1,p_1,q_1$ by $\theta, p,q$ without loss of generality we can assume that $\varphi = \varepsilon_1\theta^p$, $\psi = \varepsilon_2\theta^q$, where $|\theta|=\infty$, and $\gcd(p,q)=1$. Equalities (\ref{main equalities}) can be rewritten in the following form
\begin{align*}
\varepsilon_1^{\varphi_{11}-1}\varepsilon_2^{\varphi_{21}}\theta^{p(\varphi_{11}-1)+q\varphi_{21}} &= E,\\
\varepsilon_1^{\varphi_{12}}\varepsilon_2^{\varphi_{22}-1}\theta^{p\varphi_{12}+q(\varphi_{22}-1)} &= E,\\
\varepsilon_1^{\psi_{11}-1}\varepsilon_2^{\psi_{21}}\theta^{p(\psi_{11}-1)+q\psi_{21}} &= E,\\
\varepsilon_1^{\psi_{12}}\varepsilon_2^{\psi_{22}-1}\theta^{p\psi_{12}+q(\psi_{22}-1)} &= E.
\end{align*}
Since $|\theta|=\infty$, there is no power of $\theta$ which is equal to $-E$. Therefore from the last system of equalities we have the system
$$
\begin{cases}
p(\varphi_{11}-1)+q\varphi_{21}=0,\\
p\varphi_{12}+q(\varphi_{22}-1)=0,\\
p(\psi_{11}-1)+q\psi_{21}=0,\\
p\psi_{12}+q(\psi_{22}-1)=0.
\end{cases}
$$
From this system we have $\varphi_{11} = 1 - \dfrac{q}{p}\varphi_{21}$, $\varphi_{22} = 1 - \dfrac{p}{q}\varphi_{12}$ and
$$\varphi=\begin{pmatrix}1 - \dfrac{q}{p}\varphi_{21}&\varphi_{12}\\
\varphi_{21}&1 - \dfrac{p}{q}\varphi_{12}
\end{pmatrix}.$$

If ${\rm det}(\varphi)=-1$, then $\dfrac{p}{q}\varphi_{12}+\dfrac{q}{p}\varphi_{21}=2$. Hence 
$${\rm tr}(\varphi)=\varphi_{11}+\varphi_{22}=1 - \dfrac{q}{p}\varphi_{21}+1 - \dfrac{p}{q}\varphi_{12}=0.$$
Since ${\rm det}(\varphi)=-1$ and ${\rm tr}(\varphi)=0$ from Lemma~\ref{orders} it follows that $|\varphi|=2$ what contradicts the equality $|\theta|=\infty$.

If ${\rm det}(\varphi)=1$, then $\dfrac{p}{q}\varphi_{12} + \dfrac{q}{p}\varphi_{21} = 0$, therefore $\varphi_{12} = -\dfrac{q^2}{p^2}\varphi_{21}$ and
$$\varphi=\begin{pmatrix}
1 - \dfrac{q}{p}\varphi_{21} & -\dfrac{q^2}{p^2}\varphi_{21} \\
\varphi_{21} & 1 + \dfrac{q}{p}\varphi_{21}
\end{pmatrix}.$$
Since $\varphi$ is an integral matrix, $\varphi_{21}=sp^2$ for appropriate $s\in \mathbb{Z}$. Hence
$$\varphi=\begin{pmatrix}
1 - spq & -sq^2 \\
sp^2 & 1 + spq
\end{pmatrix}=
\begin{pmatrix}
1 - pq & -q^2 \\
p^2 & 1 + pq
\end{pmatrix}^{s}.$$
Since ${\rm tr}(\varphi)=2$, the matrix $\varphi$ is parabolic, and it has only a single eigenvalue, which is either $1$ or $-1$. Since $\varphi=\varepsilon_1\theta^p$, one of the following conditions hold.
\begin{itemize}
\item If $\varepsilon_1=-1$, then ${\rm tr}(\theta)=-2$ and $p$ is odd.
\item If $\varepsilon_1=1$, then ${\rm tr}(\theta)=\pm2$ and $p$ is even.
\end{itemize}
If necessary, replacing the matrix $\theta$ by the matrix $-\theta$, we can assume that $\varepsilon_1=\varepsilon_2=1$, and $\varphi=\theta^p, \psi=\theta^q$. Due to the fact that $p,q$ are coprime, and properties of parabolic matrices, we have $s = mp$ and
$$\theta =
\begin{pmatrix}
1 - pq & -q^2 \\
p^2 & 1 + pq
\end{pmatrix}.
$$
Hence
\begin{align*}
\varphi =
\begin{pmatrix}
1 - mp^2q & -mpq^2 \\
mp^3 & 1 + mp^2q
\end{pmatrix},&&\psi =
\begin{pmatrix}
1 - mpq^2 & -mq^3 \\
mp^2q & 1 + mpq^2
\end{pmatrix}.
\end{align*}
This is exactly case 1.2 from Table 1. The case $|\varphi|=\infty$ is finished. \hfill$\square$
 
Theorem~\ref{ttthhhtt}, Lemma~\ref{th3} and Lemma~\ref{conclemma} imply the following main result of the paper.

\medskip\noindent \textbf{Theorem A.} {\it Let $(A,\oplus,\odot)$ be a $\lambda$-homomorphic brace with $A_{\oplus}=\mathbb{Z}^2$. Then the operations in this brace are given by formulas
\begin{align*}\begin{pmatrix}a_1\\a_2\end{pmatrix}\oplus\begin{pmatrix}b_1\\b_2\end{pmatrix}=\begin{pmatrix}a_1+b_1\\a_2+b_2\end{pmatrix},&&\begin{pmatrix}a_1\\a_2\end{pmatrix}\odot\begin{pmatrix}b_1\\b_2\end{pmatrix}=\begin{pmatrix}a_1\\a_2\end{pmatrix}+\varphi^{a_1}\psi^{a_2}\begin{pmatrix}b_1\\b_2\end{pmatrix},
\end{align*}
where $\varphi,\psi\in{\rm GL}_2(\mathbb{Z})$ are matrices which satisfy one of the conditions listed in Table 1.}

\medskip\noindent \textbf{Proof.} Let $A_{\oplus}=\mathbb{Z}^2=\langle x,y~|~[x,y]=1\rangle$ be the free abelian group with the generators $x,y$. Since $A$ is $\lambda$-homomorphic, for arbitrary $a_1,a_2\in \mathbb{Z}$ we have $\lambda_{a_1x+a_2y}=\varphi^{a_1}\psi^{a_2}$, where $\varphi=\lambda_x, \psi=\lambda_y$. Hence, the operation $\odot$ is given by the formula
$$\begin{pmatrix}a_1\\a_2\end{pmatrix}\odot\begin{pmatrix}b_1\\b_2\end{pmatrix}=\begin{pmatrix}a_1\\a_2\end{pmatrix}+\varphi^{a_1}\psi^{a_2}\begin{pmatrix}b_1\\b_2\end{pmatrix}.
$$
From Lemma~\ref{th3} and Lemma~\ref{conclemma} it follows that $\varphi\psi=\psi\varphi$ and  conditions (\ref{main equalities}) hold. By Theorem~\ref{ttthhhtt} all pairs of such matrices are described in Table 1.\hfill$\square$

{\footnotesize
\begin{spacing}{0.5}

\end{spacing}}

~\\
Timur Nasybullov (timur.nasybullov@mail.ru), Igor Novikov (i.novikov1@g.nsu.ru)\\
Novosibirsk State University, Pirogova 1, 630090 Novosibirsk, Russia

\end{document}